\documentclass{amsart}[11pt,amssymb]
\DeclareMathSymbol{\twoheadrightarrow} {\mathrel}{AMSa}{"10}

\def\Q{{\mathbf Q}}
        \def\PP{\mathcal{P}}

        \def\CC{\mathfrak{C}}
        \def\CV{\mathcal{C}}
\def\Z{{\mathbf Z}}
\def\C{{\mathbf C}}
\def\R{{\mathbf R}}
\def\F{{\mathbf F}}

\def\H{{\mathrm H}}
\def\ST{{\mathbf{S}}}
\def\A{{\mathbf A}}
\def\Sn{{\mathbf S}_n}
\def\An{{\mathbf A}_n}

\def\Gal{\mathrm{Gal}}

\def\Perm{\mathrm{Perm}}

\def\Lie{\mathrm{Lie}}

\def\End{\mathrm{End}}

\def\Aut{\mathrm{Aut}}
\def\Hom{\mathrm{Hom}}
\def\Mat{\mathrm{Mat}}
\def\cl{\mathrm{cl}}

    \def\RR{\mathfrak{R}}
\def\I{\mathrm{Id}}
\def\J{{\mathcal J}}
     \def\B{\mathcal{B}}

    \def\W{{\mathcal W}}

        \def\K_a{\bar{K}}

\def\dim{\mathrm{dim}}
\def\O{\mathfrak{O}}
\def\P{{\mathbf P}}

\def\X{{\mathcal X}}

\newtheorem{thm}{Theorem}[section]
\newtheorem{lem}[thm]{Lemma}
\newtheorem{cor}[thm]{Corollary}

\theoremstyle{definition}
\newtheorem{defn}[thm]{Definition}
\newtheorem{ex}[thm]{Example}
\newtheorem{exs}[thm]{Examples}

\newtheorem{rem}[thm]{Remark}
\newtheorem{rems}[thm]{Remarks}

\title[Superelliptic jacobians]{Superelliptic jacobians}
\author{Yuri G. Zarhin}
\address{Department of Mathematics, Pennsylvania
State University, University Park, PA 16802, USA}
\email{zarhin\char`\@math.psu.edu}

\begin{document}
\maketitle

\section{Introduction}

Throughout this paper $K$ is a field of characteristic zero,
$\K_a$ its algebraic closure and $\Gal(K)=\Aut(\K_a/K)$ the
absolute Galois group of $K$.

We write $k$ for a field of characteristic zero and $k(t)$ for the
field of rational functions over $k$ in independent variable $t$.

 If $X$ is an abelian variety over
$\bar{K}$ then we write $\End(X)$ for the ring of all its
$\bar{K}$-endomorphisms and $\End^{0}(X)$ for the corresponding
$\Q$-algebra $\End(X)\otimes\Q$; the notation $1_X$ stands for the
identity automorphism of $X$. If $m$ is a positive integer then we
write $X_m$ for the kernel of multiplication by $m$ in
$X(\bar{K})$. It is well-known \cite{Mumford} that  $X_m$ is a
free $\Z/m\Z$-module of rank $2\dim(X)$. If $X$ is defined over
$K$ then $X_m$ is a Galois submodule in $X(\bar{K})$.

\begin{defn}
Suppose that $K$ contains a field $k$ and $X$ is an abelian
variety defined over $K$. We say that $X$ is isotrivial over $k$
if there exists an abelian variety $X_0$ over $\bar{k}$ such that
$X$ and $X_0$ are isomorphic over $\bar{K}$.

We say that $X$ is {\sl completely non-isotrivial} over $k$ if for
every abelian variety $W$ over $\bar{k}$ there are no non-zero
homomorphisms between $X$ and $W$ over $\bar{K}$.
\end{defn}

Let $f(x)\in K[x]$ be a polynomial of degree $n\ge 3$ with
coefficients in $K$ and without multiple roots, $\RR_f\subset
\bar{K}$ the ($n$-element) set of roots of $f$ and
$K(\RR_f)\subset \bar{K}$ the splitting field of $f$. We write
$\Gal(f)=\Gal(f/K)$ for the Galois group $\Gal(K(\RR_f)/K)$ of
$f$; it permutes roots of $f$ and may be viewed as a certain
permutation group of $\RR_f$, i.e., as as a subgroup of the group
$\Perm(\RR_f)\cong\Sn$ of permutation of $\RR_f$. ($\Gal(f)$ is
transitive if and only if $f$ is irreducible.)

 Suppose that  $p$ is a prime that does not divide $n$ and
 a positive integer $q=p^r$ is a power of $p$
 then we write $C_{f,q}$
for the superelliptic $K$-curve $y^q=f(x)$ and $J(C_{f,q})$ for
its jacobian. Clearly, $J(C_{f,q})$ is an abelian variety that is
defined over $K$ and
$$\dim(J(C_{f,q}))=\frac{(n-1)(q-1)}{2}.$$
In a series of papers \cite{ZarhinMRL,ZarhinMMJ,ZarhinBSMF},
\cite{ZarhinCrelle,ZarhinCamb,ZarhinSb,ZarhinM}, the author
discussed the structure of $\End^0(J(C_{f,q}))$, assuming that
$n\ge 5$ and $\Gal(f)$ is, at least, doubly transitive. In
particular, he proved that if $n \ge 5$ and $\Gal(f)$ coincides
either with  full symmetric group $\Sn$ or with  alternating group
$\An$ then $\End^0(J(C_{f,q}))$ is (canonically) isomorphic to a
product $\prod_{i=1}^r\Q(\zeta_{p^i})$ of cyclotomic fields. (If
$q=p$ then we proved that $\End(J(C_{f,q}))=\Z[\zeta_p]$.)

In this paper we discuss the remaining cases when $n=3$ or $4$ and
$\Gal(f)$ is a doubly transitive subgroup of $\Sn$, i.e., the
cases when
$$n=3,\ \Gal(f)=\ST_3, \ \dim(J(C_{f,q}))=q-1$$
 and
 $$n=4,\ \Gal(f)=\ST_4 \text{ or } \A_4,\  \dim(J(C_{f,q}))=\frac{3(q-1)}{2}.$$
In those cases if $q=p$ then we have
 $$1_{J(C_{f,p})}\in \Q(\zeta_p)\subset \End^0(J(C_{f,p})).$$
 If $n=3$ this means that the $(p-1)$-dimensional abelian variety $J(C_{f,p})$ admits
 multiplication by the field $\Q(\zeta_p)$ of degree $p-1$. One may prove (see Sect. \ref{abvar})
  that either $J(C_{f,p})$ is an abelian variety of
 CM-type or $\End^0(J(C_{f,p}))=\Q(\zeta_p)$ (and even $\End(J(C_{f,p}))=\Z[\zeta_p]$).

If $n=4$ this means that the $3(p-1)/2$-dimensional  abelian
variety $J(C_{f,p})$ admits multiplication by the field
$\Q(\zeta_p)$ of degree $p-1$. One may prove (see Sect.
\ref{abvar}) that either $J(C_{f,p})$ contains an abelian
subvariety of CM-type (of positive dimension) or
$\End^0(J(C_{f,p}))=\Q(\zeta_p)$ (and even
$\End(J(C_{f,p}))=\Z[\zeta_p]$).

\begin{thm}
\label{primep}
 Suppose that $k$ is algebraically closed and $K=k(t)$. Assume that $(n,\Gal(f))$ enjoy one of
the following two properties:
\begin{itemize}
\item $n=3$ and $\Gal(f)=\ST_3$. \item $n=4, \Gal(f)=\ST_4 \text{
or } \A_4$.
\end{itemize}
Suppose that $p$ is an odd prime that does not divide $n$.  Then
$J(C_{f,p})$ is completely non-isotrivial and
$\End(J(C_{f,p}))=\Z[\zeta_p]$.
\end{thm}

\begin{thm}
\label{oddp} Suppose that $k$ is algebraically closed and
$K=k(t)$. Assume that $(n,\Gal(f))$ enjoy one of the following two
properties:
\begin{itemize}
\item $n=3$ and $\Gal(f)=\ST_3$. \item $n=4, \Gal(f)=\ST_4 \text{
or } \A_4$.
\end{itemize}
Suppose that $p$ is a prime that does not divide $n$. Let $r$ be a
positive integer and $q=p^r$. If $p$ is odd  then  $J(C_{f,q})$ is
completely non-isotrivial and $\End^0(J(C_{f,q}))$ is
(canonically) isomorphic to a product
$\prod_{i=1}^r\Q(\zeta_{p^i})$ of cyclotomic fields.
\end{thm}

\begin{thm}
\label{deg3} Suppose that $k$ is algebraically closed and
$K=k(t)$. Assume that $n=3$ and $\Gal(f)=\ST_3$.

Then:

\begin{itemize}
\item
 $\End^0(J(C_{f,4}))=\Q\times \Mat_2(\Q(\sqrt{-1}))$. In addition $J(C_{f,4})$ is
isogenous to a product of the elliptic curve $C_{f,2}:y^2=f(x)$
and the square of the elliptic curve $y^3=x^3-x$.
 \item Let $r>2$ be an integer and $q=2^r$.
  Then
  $\End^0(J(C_{f,q}))$ is (canonically) isomorphic
to a product $\Q\times
\Mat_2(\Q(\sqrt{-1}))\times\prod_{i=3}^r\Q(\zeta_{p^i})$. In
addition, $J(C_{f,q})$ is isogenous to a product of  a completely
non-isotrivial abelian variety, the elliptic curve
$C_{f,2}:y^2=f(x)$ and the square of the elliptic curve
$y^3=x^3-x$.
\end{itemize}
\end{thm}

In the case of arbitrary $n\ge 3$ and doubly transitive $\Gal(f)$
we have the following non-isotriviality assertion.
\begin{thm}
\label{double}
 Suppose that $k$ is algebraically closed, $K=k(t)$, a prime $p$ does not divide $n$ and
 $\Gal(f)$ acts doubly transitively on $\RR_f$. Let $r$ be a positive integer and $q=p^r$.
  Let us assume additionally that $(n,p)\ne
(3,2)$.

 Then $J(C_{f,q})$ is completely non-isotrivial.
\end{thm}


The paper is organized as follows. In Section \ref{ellcurves} we
discuss elliptic curves, using elementary means. Section
\ref{abvar} contains auxiliary results concerning abelian varietis
with big endomorphism fields. Section \ref{complab} deals with
complex abelian varieties. In Section \ref{superjac} we discuss in
detail superelliptic curves and their jacobians. It also contains
the proof of main results.

I am grateful to J.-L. Colliot-Th\'el\`ene and to the referee for
helpful comments. My special thanks go to Dr. Boris Veytsman for
his help with \TeX nical problems.

\section{Elliptic curves}
\label{ellcurves}
 We start with the case of $n=3$ and an elliptic
curve $C_{f,2}:y^2=f(x)$ where $f(x)$ is a cubic polynomial
without multiple roots. Then $J(C_{f,2})=C_{f,2}$ is a
one-dimensional abelian variety. If
$\{\alpha_1,\alpha_2,\alpha_3\}\subset \bar{K}$ is the set of
roots of $f(x)$ then the group of points of order $2$ on $C_{f,2}$
is $\{\infty, (\alpha_1,0), (\alpha_2,0), (\alpha_3,0)\}$. Here
$\infty\in C_{f,2}(K)$ is the zero of group law on $C_{f,2}$.

Recall that the only doubly transitive subgroup of $\ST_3$ is
$\ST_3$ itself.

\begin{exs}
Let us put $K=\Q$.
\begin{enumerate}
\item
 The polynomial $f(x)=x^3-x-1$ has Galois group $\ST_3$. The
 corresponding
elliptic curve $C_{f,2}:y^2=x^3-x-1$ has non-integral
$j$-invariant $1728\cdot \frac{-4}{23}$ and therefore
$\End(J(C_{f,2}))=\Z$. \item The polynomial
 $x^3-2$ has Galois group $\ST_3$. The
elliptic curve $y^2=x^3-2$  admits an automorphism of
multiplicative order $3$ and
$\End(J(C_{f,2}))=\Z[\frac{-1+\sqrt{-3}}{2}]$.
\end{enumerate}
\end{exs}

\begin{ex}(See \cite[Sect. 4]{EZ}.)
Let $p>3$ be a prime that is congruent to $3$ modulo $8$,
$$\omega:=\frac{-1+\sqrt{-p}}{2}, \ \alpha: = j(\omega)\in \C, \ K:=\Q(j(\omega))\subset \C.$$
Let us consider the cubic polynomial $$h_p(x):=x^3 -
\dfrac{27\alpha}{4(\alpha-1728)}x -
\dfrac{27\alpha}{4(\alpha-1728)}\in K[x].
$$ without multiple roots.
Then the $j$-invariant of the elliptic curve
$C_{h_p,2}:y^2=h_p(x)$ coincides with $\alpha$ and therefore the
endomorphism ring of $J(C_{h_p,2})$ coincides with
$\Z[\frac{-1+\sqrt{-p}}{2}]$. However, the Galois group of
$h_p(x)$ over $K$ is $\ST_3$ \cite[Sect. 4, Lemma 4.6]{EZ}.
\end{ex}

\begin{ex}
Suppose that $K=\C(t)$. The polynomial $x^3-x-t$ has Galois group
$\ST_3$. The elliptic curve $y^2=x^3-x-t$ has transcendental
$j$-invariant $\frac{-4\cdot 1728}{27t^2-4}$ and therefore its
endomorphism ring is $\Z$.
\end{ex}

\begin{lem}
\label{j}
 Suppose that  $k$ is an algebraically closed field of characteristic zero and $K\supset k$.
 Let $f(x)\in K[x]$ be a cubic
polynomial with Galois group $\ST_3$. Then the $j$-invariant of
the elliptic curve $C_{f,2}$ does not lie in $k$, i.e., $C_{f,2}$
is not isotrivial.
\end{lem}

\begin{proof}
Indeed, suppose that there exists an elliptic curve $E$ over $k$
such that $E$ and $C_{f,2}$ are isomorphic over $\bar{K}$, i.e.,
$C_{f,2}$ is a twist of $E$. Clearly, all $\bar{K}$-automorphisms
of $E$ are defined over $k$ and $\Aut(E)$ is a cyclic group of
order $2,4$ or $6$. It follows that the first Galois cohomology
group $H^1(K,\Aut(E))$ coincides with the group of continuous
homomorphisms of $\Gal(K)$ to $\Aut(E)$. It follows that there is
a cyclic extension $L/K$ and an embedding
$c:\Gal(L/K)\hookrightarrow \Aut(E)$ such that $C_{f,2}$ is the
$L/K$-twist of $E$ attached to $c$. Notice that all torsion points
of $E$ are defined over $k$ and therefore lie in $E(K)$. This
implies that the Galois action on the points of order $2$ of
$C_{f,2}$ factors through $\Gal(L/K)$; in particular, the
corresponding Galois image is commutative. Since $\ST_3$ is
non-commutative, we get a desired contradiction.
\end{proof}


\section{Abelian varieties}
\label{abvar} Let $X$ be an abelian variety of positive dimension
over an algebraically closed field $\bar{K}$ of characteristic
zero, let $E$ be a number field and let $i:E
\hookrightarrow\End^0(X)$ be an embedding such that $i(1)=1_X$. It
is well-known that $[E:\Q]$ divides $2\dim(X)$
\cite{Mumford,Shimura}; if $2\dim(X)=[E:\Q]$ then $X$ is an
abelian variety of CM-type.

 We write $\End^0(X,i)$ for the centralizer of $i(E)$
in $\End^0(X)$; it is known \cite{ZarhinL} that $\End^0(X,i)$ is a
finite-dimensional semisimple $E$-algebra.

\begin{thm}
\label{bigend} If $\dim(X)=[E:\Q]$ then either $X$ is an abelian
variety of CM-type or $\End^0(X,i)=i(E)\cong E$.

If $2\dim(X)=3[E:\Q]$ then either $X$ contains a non-zero abelian
subvariety of CM-type or $\End^0(X,i)=i(E)\cong E$.
\end{thm}

\begin{proof}
Recall that if $Z$ is an abelian variety of positive dimension $d$
then every {\sl semisimple commutative} $\Q$-subalgebra of
$\End^0(Z)$ has $\Q$-dimension $\le 2d$. The algebra  $\End^0(Z)$
contains a $2d$-dimensional semisimple commutative $\Q$-subalgebra
if and only if $Z$ is an abelian variety of CM-type.

By theorem 3.1 of \cite{ZarhinM},
$$\dim_E(\End^0(X,i))\le \left(\frac{2\dim(X)}{[E:\Q]}\right)^2;$$
if the equality holds then $X$ is an abelian variety of CM-type.
If $\dim_E(\End^0(X,i))=1$ then $\End^0(X,i)=i(E)\cong E$.

 So, in the course of the proof we may assume that either
$\dim(X)=[E:\Q]$ and $1<\dim_E(\End^0(X,i))<4$ or
$2\dim(X)=3[E:\Q]$ and $1<\dim_E(\End^0(X,i))<9$.

Now assume that $\dim(X)=[E:\Q]$. Then $\dim_E(\End^0(X,i))$=$2$
or $3$. Since every semisimple algebra of dimension $2$ or $3$
over a field is commutative, $\End^0(X,i)$ is commutative  and its
$\Q$-dimension is either $2[E:\Q]=2\dim(X)$ or
$$3[E:\Q]=3\dim(X)>2\dim(X).$$ It follows that
$\dim_{\Q}\End^0(X,i)=2\dim(X)$ and therefore $X$ is an abelian
variety of CM-type. This proves the first assertion of the
Theorem.

Now assume that $2\dim(X)=3[E:\Q]$ and $1<\dim_E(\End^0(X,i))<9$.
So, $$2\le\dim_E(\End^0(X,i))\le 8.$$ If $\dim_E(\End^0(X,i))=3$
then $\End^0(X,i)$ is commutative  and its $\Q$-dimension is
$3[E:\Q]=2\dim(X)$ and therefore $X$ is an abelian variety of
CM-type.

So, further we assume that $\dim_E(\End^0(X,i))\ne 3$. Clearly,
there exists a positive integer $a$ such that
$$\dim(X)=3a, \ 2\dim(X)=6a, \ [E:\Q]=2a.$$
First, suppose that $\End^0(X,i)$ is not simple. Then $X$ is
isogenous to a product $\prod_{s\in S}X_s$ of abelian varieties
$X_s$ provided with embeddings
$$i_s:E \hookrightarrow \End^0(X_s), \ i_s(1)=1_{X_s};$$
in addition, the  $E$-algebra $\End^0(X_s)$ is simple for all $s$
and the $E$-algebras $\End^0(X,i)$ and $\prod_{s\in S}\End^0(X_s)$
are isomorphic \cite[Remark 1.4]{ZarhinL}. It follows that
$$6a=2\dim(X)=2\sum_{s\in S}\dim(X_s)$$ and $2a=[E:\Q]$ divides
$2\dim(X_s)$ for all $s$. We conclude that $\dim(X_s)=a$ or $2a$.
It follows that $S$ is a either two-element set  $\{s_1,s_2\}$
with
$$\dim(X_{s_1})=a, \ \dim(X_{s_2})=2a$$
or a three-element set  $\{s_1,s_2,s_3\}$ with
$$\dim(X_{s_1})= \dim(X_{s_2})=\dim(X_{s_3})=a.$$
In both cases $2\dim(X_{s_1})=2a=[E:\Q]$ and therefore $X_{s_1}$
is an abelian variety of CM-type.

So, further we may assume that $\End^0(X,i)$ is a {\sl simple}
$E$-algebra.

 If $\dim_E(\End^0(X,i))=2$ then $\End^0(X,i)$ is
commutative and therefore is a number field of degree $4a$ over
$\Q$. Since $4a$ does not divide $6a=2\dim(X)$, we get a
contradiction.

 If $\dim_E(\End^0(X,i))>3$ then
 $\dim_E(\End^0(X,i))>2\dim(X)$ and
therefore $\End^0(X,i)$ is noncommutative. Since $\End^0(X,i)$ is
simple,  $\dim_E(\End^0(X,i))$ is {\sl not} square-free. It
follows that $\dim_E(\End^0(X,i))$ is either $4$ or $8$ and
therefore
$$\dim_{\Q}(\End^0(X,i))=[E:\Q]\cdot\dim_E(\End^0(X,i))$$ is
either $8a$ or $16a$ respectively. Since none of them divides
$6a=2\dim(X)$, the simple $E$-algebra $\End^0(X,i)$ is not a
division algebra. (Here we use that $\bar{K}$ has characteristic
zero.) It follows that there exist an integer $r>1$, a division
algebra $D$ over $E$ and an isomorphism $\Mat_r(D)\cong
\End^0(X,i)$ of $E$-algebras. It follows that
$$r^2\dim_E(D)=\dim_E(\Mat_r(D))=\dim_E(\End^0(X,i))$$
is either $4$ or $8$. Clearly, there exist an abelian variety $Y$
of dimension $\frac{1}{2}\dim(X)=\frac{3}{2}a$, an isogeny $Y^2\to
X$ and an embedding
$$\iota:D \hookrightarrow \End^0(Y), \ \iota(1)=1_Y.$$
It follows that $1_Y\in \iota(E)\subset \End^0(Y)$ and the degree
$2a$ of the number field $\iota(E)\cong E$ must divide
$2\dim(Y)=3a$. In other words, $2a$ divides $3a$. Contradiction.
\end{proof}

\begin{rem}
The case when $E$ is a (product of) totally real field(s) with
$\dim_{\Q}(E)=\dim(X)$ was earlier discussed in \cite[pp.
557--558]{Ribet2}.
\end{rem}

\begin{rem}
\label{hom} Suppose that $X$ and $Y$ are absolutely simple abelian
varieties over $\bar{K}$ such that the $\Q$-algebras $\End^0(X)$
and $\End^0(Y)$ are {\sl not} isomorphic. Then $X$ and $Y$ are not
isogenous and the absolute simplicity implies that there are no
non-zero homomorphisms between $X$ and $Y$. It follows that
$\End^0(X\times Y)=\End^0(X)\times \End^0(Y)$.
\end{rem}

\begin{lem}
\label{nonconst}
 Suppose that $k$ is a field of characteristic zero and
$K=k(t)$. Let $X$ be an abelian variety over $K$, let $B \subset
\A^1_k$ be a finite set of closed points on the affine $k$-line
 and let $f:\X \to \A^1_k\setminus B$ be an abelian scheme, whose
 generic fiber coincides with $X$.
For each $a\in k\setminus B\subset k=\A^1_k(k)$ we write $X_a$ for
the corresponding fiber that is an abelian variety over $k$.

Suppose that there exists a non-zero abelian variety $W$ over
$\bar{k}$ such that there exists a non-zero homomorphism from $W$
to $X$ that is defined over $\bar{K}$. Then for any pair of
distinct points $a_1,a_2\in k\setminus B$ there exists a non-zero
 homomorphism between  abelian varieties $X_{a_1}$ and $X_{a_2}$
that is defined over $\bar{k}$.
\end{lem}

\begin{proof}
Recall that one may view closed points of $\A^1_k$ as Galois
orbits in $\bar{k}$. Replacing $k$ by $\bar{k}$, $K$ by
$\bar{k}(t)$, $X$ by $X\times_{k(t)}\bar{k}(t)$, $B$ by the union
$\bar{B}$ of all Galois orbits from $B$ and $\X \to
\A^1_k\setminus B$ by the fiber product $X \times_
{\A^1_k\setminus B}\A^1_{\bar{k}}\setminus \bar{B}\to
\A^1_{\bar{k}}\setminus \bar{B}$, we may and will assume that $k$
is algebraically closed, i.e., $k=\bar{k}$. We may assume that $W$
is simple and a non-zero homomorphism $W \to X$ is an embedding
that is defined over $\bar{K}$.

Notice that both $W$ and $X$ are defined over $K$. In addition,
all torsion points of "constant" $W$ are defined over $K$ (in
fact, they are defined even over $k$). Let $L=K(X_3)$ be the field
of definition of all points of order $3$ on $X$; clearly, $L$ is a
finite Galois extension of $K=k(t)$. It follows from results of
\cite{Silverberg} that all homomorphisms between $W$ and $X$ are
defined over $L$. In particular, there is a closed embedding
$i:W\hookrightarrow X$ that is defined over $L$.

It follows from N\'eron-Ogg-Shafarevich criterion \cite{SerreTate}
that $L/k(t)$ is unramified outside $\infty$ and $B$.

This means that if $U$ is the normalization of $\A^1_k\setminus B$
in $L$ then the natural map $\pi:U\to\A^1_k\setminus B$ is
\'etale; clearly, it is a surjective map between the sets of
$k$-points. The field $k(U)$ of rational functions on $U$
coincides with $L$.

Let $\X_U\to U$ be the pullback of $\X$ with respect to $\pi$; it
is an abelian scheme over $U$, whose generic fiber coincides with
$X\times_K L$. If $u\in U(k)$ and $$a=\pi(u)\in (\A^1_k\setminus
B)(k)=k\setminus B$$ then the corresponding closed fiber $X_{u}$
of $X_U$ coincides with $X_s$.

Recall that $X$ contains an abelian $L$-subvariety isomorphic to
$W$. Let $\W$ be the schematic closure of $W$ in $\X_U$ \cite[p.
55]{Neron}. It follows from \cite[Sect. 7.1, p. 175, Cor.
6]{Neron} that $\W$ is a N\'eron model of $W$ over $U$. On the
other hand, since $W$ is "constant", its N\'eron model over $U$ is
isomorphic to $W\times U$. This implies that the group $U$-schemes
$\W$ and $W\times U$ are isomorphic. Since $W\times U$ is a closed
$U$-subscheme of $X_U$,  the corresponding closed fiber $X_{u}$ of
$X_U$ contains an abelian subvariety  isomorphic to $W$ for all
$k$-point $u\in U(k)$. It follows that $X_s$ contains an abelian
subvariety  isomorphic to $W$ for all $a \in k\setminus B$.

Now let $a_1,a_2$ be two distinct elements of $k\setminus B$. It
follows from Poincar\'e's complete reducibility theorem
\cite[Sect. 18, Th. 1 on p. 173]{Mumford} that there exists a
surjective homomorphism $X_{a_1}\twoheadrightarrow W$. Now the
composition
$$X_{a_1}\twoheadrightarrow W\hookrightarrow X_{a_2}$$
is a non-zero homomorphism between $X_{a_1}$ and $X_{a_2}$.
\end{proof}

\section{Complex abelian varieties}
\label{complab} \label{MT} Let $Z$ be a complex abelian variety of
positive dimension and  $\CC_Z$  the center of the semisimple
finite-dimensional $\Q$-algebra $\End^0(Z)$. We write
$\Omega^1(Z)$ for the space of differentials of the first kind on
$Z$; it is a $\dim(Z)$-dimensional $\C$-vector space.

Let $E$ be a subfield of $\End^0(Z)$ that contains $1_Z$. The
degree $[E:\Q]$ divides $2\dim(Z)$. We write
$$\iota: E \subset \End^0(Z)$$
for the inclusion map and $\Sigma_E$ for the set of all field
embeddings $\sigma: E \hookrightarrow \C$. The centralizer
$\End^0(Z,\iota)$ is a semisimple $E$-algebra.

  It is well-known that
$$\C_{\sigma}:=E\otimes_{E,\sigma}\C=\C, \quad E_{\C}=E\otimes_{\Q}\C=\prod_{\sigma\in \Sigma_E}
E\otimes_{E,\sigma}\C=\prod_{\sigma\in \Sigma_E}\C_{\sigma}.$$ By
functoriality, $\End^0(Z)$ and therefore $E$ act on $\Omega^1(Z)$
and therefore provide $\Omega^1(Z)$ with a natural structure of
$E\otimes_{\Q}\C$-module \cite[p. 341]{ZarhinM}. Clearly,
$$\Omega^1(Z)=\bigoplus_{\sigma\in
\Sigma_E}\C_{\sigma}\Omega^1(Z)=\oplus_{\sigma\in
\Sigma_E}\Omega^1(Z)_{\sigma}$$ where
$\Omega^1(Z)_{\sigma}:=\C_{\sigma}\Omega^1(Z)=\{x \in
\Omega^1(Z)\mid ex=\sigma(e)x \quad \forall e\in E\}$. Let us put
$n_{\sigma}=n_{\sigma}(Z,E)=\dim_{\C_{\sigma}}\Omega^1(Z)_{\sigma}=\dim_{\C}\Omega^1(Z)_{\sigma}$.
  Let us put
$$n_{\sigma}:=\dim_{\C}(\Omega^1(Z)_{\sigma})=\dim_{\C}(\Omega^1(Z)_{\sigma}) \eqno (1).$$

\begin{thm}
\label{mult} If  $E/\Q$ is Galois, $E\supset\CC_Z$ and $\CC_{Z}
\ne E$ then there exists a nontrivial automorphism $\kappa: E \to
E$ such that $n_{\sigma}=n_{\sigma\kappa}$ for all $\sigma\in
\Sigma_E$.
\end{thm}

\begin{proof}
See \cite[Th. 2.3]{ZarhinCamb}.
\end{proof}

\begin{thm}
\label{cyclmult} Suppose that there exist a prime $p$, a positive
integer $r$,  the prime  power $q=p^r$ and an integer $n\ge 3$
enjoying the following properties:
\begin{itemize}
 \item[(i)] $E=\Q(\zeta_q)\subset\C$ where $\zeta_q\in\C$ is
a primitive $q$th root of unity; \item[(ii)] $n$ is not divisible
by $p$, i.e. $n$ and $q$ are relatively prime: \item[(iii)] Let
$i<q$ be a positive integer that is not divisible by $p$ and
$\sigma_i:E=\Q(\zeta_q)\hookrightarrow \C$ the embedding that
sends $\zeta_q$ to $\zeta_q^{-i}$. Then
$n_{\sigma_i}=\left[\frac{ni}{q}\right]$.
\end{itemize}

Then
\begin{itemize}
\item[(a)] If $E\supset\CC_Z$ then $\CC_{Z}=\Q(\zeta_q)$.

\item[(b)] If
$$\dim_E(\End^0(Z,\iota))=
\left(\frac{2\dim(Z)}{[E:\Q]}\right)^2$$ then $n=3,p=2,q=4$.
\end{itemize}
\end{thm}

\begin{proof}
 If $n>3$ then  (a) is proven in \cite[Cor. 2.2 on p.
342]{ZarhinM}. So, in order to prove (a), let us assume that $n=3$
and $\CC_{Z}\ne\Q(\zeta_q)$.

 Clearly, $\{\sigma_i\}$ is the collection $\Sigma$ of
all embeddings $\Q(\zeta_q)\hookrightarrow \C$.  By (iii),
$n_{\sigma_i}=0$ if and only if $1\le i \le [\frac{q}{3}]$.
Suppose that $\CC_{Z} \ne \Q(\zeta_q)$.
 It follows from Theorem \ref{mult}  that
there exists a non-trivial field automorphism
 $\kappa: \Q[\zeta_q] \to \Q[\zeta_q]$ such that
 for all $\sigma\in \Sigma$ we have
$n_{\sigma}=n_{\sigma\kappa}$.
 Clearly, there exists an integer $m$ such that $p$ does {\sl not}
 divide $m$,
$1<m<q$ and $\kappa(\zeta_q)=\zeta_q^m$.

If $i$ is an integer then we write $\bar{i}\in \Z/q\Z$ for its
  residue modulo $q$.
Clearly, $n_{\sigma}=0$ if and only if $\sigma=\sigma_i$ with $1
\le i \le [\frac{q}{3}]$. Since $3$ and $q$ are relatively prime,
$[\frac{q}{3}]=[\frac{q-1}{3}]$. It follows that $n_{\sigma_i}=0$
if and only if $1\le i \le [\frac{q-1}{3}]$. Clearly, the map
$\sigma \mapsto \sigma \kappa$ permutes the set $$\{\sigma_i\mid 1
\le i \le \left[\frac{q-1}{3}\right], p \mathrm{\ does\ not \
divide}\ i\}.$$

Since $\kappa(\zeta_q)=\zeta_q^m$ and
$\sigma_i\kappa(\zeta_q)=\zeta_q^{-im}$, it follows that if
$$A:=\left\{i \in \Z\mid 1 \le i \le \left[\frac{q-1}{3}\right]<q,\
p \mathrm{\ does\ not \ divide}\ i \right\}$$ then the
multiplication by $m$ in $(\Z/q\Z)^{*}=\Gal(\Q(\zeta_q)/\Q)$
leaves invariant the set $\bar{A}:=\left\{\bar{i}\in \Z/q\Z\mid
 i\in A\right\}$.
Clearly, $A$ contains $1$ and therefore $\bar{m}=m\cdot \bar{1}\in
\bar{A}$. Since $1<m<q$,
$$m=m\cdot 1\le \left[\frac{q-1}{3}\right]
 \eqno (2).$$
 Since $m\ge 2$, we conclude that
 $$q\ge 7.$$
 Let us consider the arithmetic
progression consisting  of  $2m$ integers $$[\frac{q-1}{3}]+1,
\ldots , [\frac{q-1}{3}]+2m$$ with difference $1$. All its
elements lie between $[\frac{q-1}{3}]+1$ and
$$\left[\frac{q-1}{3}\right]+2m  \le 3\left[\frac{q-1}{3}\right]\le
3\frac{q-1}{3}\le q-1.$$ Clearly, there exist exactly two elements
of $A$ say, $mc_1$ and $mc_1+m$ that are  divisible by $m$.
Clearly, $c_1$ is a positive integer and either $c_1$ or $c_1+1$
is not divisible by $p$; we put $c=c_1$ in the former case and
$c=c_1+1$ in the latter case. However, $c$ is not divisible by $p$
and
$$\left[\frac{q-1}{3}\right]<mc\le
\left[\frac{q-1}{3}\right]+2m\le q-1 \eqno(3).$$ It follows that
$mc$ does not lie in $A$ and therefore $\overline{mc}$ does not
lie in $\bar{A}$. This implies  that $\bar{c}$ also does not lie
in $\bar{A}$ and therefore $c>\left[\frac{q-1}{3}\right]$. Using
(3), we conclude that $$(m-1)\left[\frac{q-1}{3}\right]<2m$$ and
therefore
$$\left[\frac{q-1}{3}\right]<\frac{2m}{m-1}=2+\frac{2}{m-1}.$$
If $m>2$ then $m\ge 3$ and using (2), we conclude that
$$3 \le m \le \left[\frac{q-1}{3}\right]< 2+\frac{2}{m-1}< 3$$
and therefore $3<3$, which is not true. Hence $m=2$. Again using
(2), we conclude that
$$2=m \le\left[\frac{q-1}{3}\right]< 2+\frac{2}{m-1}=4$$
and therefore $\left[\frac{q-1}{3}\right]=2$ or $3$. It follows
that $q=7$ or $11$. (Since $m=2$ is prime to $q$, we have $q\ne
8$.)
  We conclude that  either
$\bar{A}=\{\bar{1},\bar{2}\}$ or $A=\{\bar{1},\bar{2},\bar{3}\}$.
In both cases $\bar{4}=2\cdot \bar{2}=m\cdot \bar{2}$ must lie in
$\bar{A}$. Contradiction. This proves (a).

Proof of (b). It follows from  Theorem 4.2 of \cite{ZarhinL} that
there exist a complex abelian variety $W$ with
$2\dim(W)=\varphi(q)=(p-1)p^{r-1}$, an embedding
$$\epsilon: E=\Q(\zeta_q)\hookrightarrow \End^0(W),
\epsilon(1)=1_W,$$ and an isogeny $\psi: W^{n-1} \to Z$ that
commutes with the action(s) of $E$. Clearly $\psi$ induces an
isomorphism of $E\otimes_{\Q}\C$-modules
$$\psi^{*}: \Omega^1(Z)\cong \Omega^1(W^{n-1})=\Omega^1(W)^{n-1}$$
where $\Omega^1(W)^{n-1}$ is the direct sum of ${n-1}$ copies of
the $E\otimes_{\Q}\C$-module $\Omega^1(W)$. This implies that if
any $\gamma \in E$ is viewed as a $\C$-linear operator in
$\Omega^1(Z)$ then it has, at most, $\dim(W)$ distinct eigenvalues
and the multiplicity of each eigenvalue is divisible by $n-1$.
Applying this observation to $\gamma=\zeta$, we conclude that the
set
$$B:=\left\{i\in\Z\mid 1\le i<q, (i,p)=1, \left[\frac{ni}{q}\right]>0\right\}=
\left\{i\in\Z\mid \frac{q}{n} < i < q=p^r, (i,p)=1\right\}$$ has,
at most $\dim(W)$ elements. Clearly, $\dim(W)=(p-1)p^{r-1}/2$ is
an integer; in particular, $q=p^r\ne 2$.

If $p/n<(p-1)/2$ then
$$\#(B)> \varphi(q)-\frac{q}{n}=(p-1)p^{r-1}-\frac{p^{r-1}p}{n}=
\left(p-1-\frac{p}{n}\right)p^{r-1}>\frac{p-1}{2}p^{r-1}=\dim(W)$$
and therefore $\#(B)>\dim(W)$, which is not the case. Hence
$$\frac{p}{n}>\frac{p-1}{2}.$$
This implies easily that $n=3$ and therefore $p\ne 3$. It follows
that $p/3>(p-1)/2$, i.e. $1/2>p/6$ and therefore $3\ge p$. This
implies that $p=2$. Now if $r\ge 3$ then $q=2^r\ge 8$ and
$2^r/6>1$.  Then $B$ contains all $2^{r-2}$ odd integers that lie
between $2^{r-1}$ and $2^r$. In addition,
$$2^{r-1}-\frac{2^r}{3}=\frac{2^r}{6}>1$$
and therefore an odd  number $2^{r-1}-1$ is greater than $2^r/3$,
i.e., $2^{r-1}-1$ lies in $B$. This implies that
$\#(B)>2^{r-2}=\varphi(q)/2=\dim(W)$, which is not the case.
\end{proof}

\begin{lem}
\label{q4n3} Suppose that
$$p=2,q=4, E=\Q(\zeta_4)=\Q(\sqrt{-1}), n=3.$$
Assume that if $i$ is an odd positive integer that is strictly
less than $4$  and $\sigma_i:E=\Q(\zeta_4)\hookrightarrow \C$ the
embedding that sends $\zeta_q$ to $\zeta_q^{-i}$ then
$n_{\sigma_i}=\left[\frac{3i}{q}\right]$.

Then:
\begin{itemize}
\item[(i)]
 $Z$ is isomorphic to a product of two mutually isogenous
elliptic curves with complex multiplication by $\Q(\sqrt{-1})$.
\item[(ii)] $Z$ is isogenous to the square of the elliptic curve
$y^2=x^3-x$. \item[(iii)]
 $\End^0(Z)\cong\Mat_2(\Q(\sqrt{-1}))$.
 \end{itemize}
\end{lem}

\begin{proof}
Clearly,
$$n_{\sigma_1}=0, n_{\sigma_3}=2, \dim(Z)=\dim(\Omega^1(Z))=n_{\sigma_1}+n_{\sigma_3}=2.$$
This means that $E$ acts on the complex vector space
$\Omega^1(Z))$ as multiplications by scalars via
$$\sigma_3:E \hookrightarrow \C.$$
Let us put $\O:=E\bigcap \End(Z)$; clearly, it is an order in $E$.
Let $\Lie(Z)$ be the tangent space to $Z$ at the origin; one may
view $\Omega^1(Z)$ as the complex vector space, which is dual to
$\Lie(Z)$. In particular, $E$ also acts on $\Lie(Z)$ as
multiplications by scalars via $\sigma_3$. Further, we identify
$E$ with its image in $\C$ via $\sigma_3$. Clearly, the natural
map
$$\O\otimes\R\to\C$$
is an isomorphism of $\R$-algebras.

The exponential surjective holomorphic map
$$\exp:\Lie(Z)\twoheadrightarrow Z(\C)$$
is a homomorphism of $\O$-modules. Recall that  he kernel $\Gamma$
of $\exp$ is a discrete lattice of rank $2\dim(Z)=4$ in $\Lie(Z)$.
Clearly, $\Gamma$ is $\O$-stable and therefore may be viewed as an
$\O$-module  of rank $2$. We have
$$Z(\C):=\Lie(Z)/\Gamma.$$
Clearly, the natural map
$$\Gamma\otimes\R\to\Lie(Z)$$
is an isomorphism of real vector spaces.
 Since $\O$ is the order in the quadratic field, it follows
from a theorem of Z. Borevich - D. Faddeev  \cite{BF} (see also
\cite[Satz 2.3]{Schoen}) that the $\O$-module $\Gamma$ splits into
a direct sum
$$\Gamma=\Gamma_1\oplus\Gamma_2$$
of rank one $\O$-modules $\Gamma_1$ and $\Gamma_2$. Clearly
$$\Lie(Z)=\Gamma_1\otimes\R\oplus\Gamma_2\otimes\R$$
and both $\Gamma_j\otimes\R$ are $\O\otimes\R=\C$-submodules in
$\Lie(Z)$ ($j=1,2$). It follows that the complex torus
$Z(\C):=\Lie(Z)/\Gamma$ is isomorphic to a product
$(\Gamma_1\otimes\R)/\Gamma_1\times (\Gamma_2\otimes\R)/\Gamma_2$
of one-dimensional complex tori $(\Gamma_1\otimes\R)/\Gamma_1$ and
$(\Gamma_2\otimes\R)/\Gamma_2$ with complex multiplication by
$\O$. One has only to recall that every one-dimensional complex
torus is an elliptic curve and two CM-elliptic curves with multiplication by
the same imaginary quadratic field are isogenous.

\end{proof}

\section{Superelliptic jacobians}
\label{superjac}
Throughout this Section $p$ is a prime, $r$ a positive integer and $q=p^r$.
We write ${\mathcal P}_q(t)$ for the polynomial $1+t+\dots + t^{q-1}\in \Z[t]$ and $\Phi_q(t)$ for the
$q$th cyclotomic polynomial $\Phi_q(t)=\sum_{i=0}^{p-1}t^{i p^{r-1}}\in \Z[t]$ of degree
 $\varphi(q)=q-p^{r-1}$. Notice that
 $$t^q \Phi_q(1/t)-\Phi_q(t)=t^q-1 \eqno (4).$$

 Let $K$ be a field of characteristic $0$ and $f(x) \in
K[x]$ be a separable polynomial of degree $n \ge 3$. We write
$\RR_f$ for the set of its roots and denote by $K(\RR_f)\subset
\bar{K}$ the corresponding splitting field. As usual, the Galois
group $\Gal(K(\RR_f)/K)$ is called the Galois group of $f$ and
denoted by $\Gal(f)$. Clearly, $\Gal(f)$ permutes elements of
$\RR_f$ and the natural map of $\Gal(f)$ into the group
$\Perm(\RR_f)$ of all permutations of $\RR_f$ is an embedding. We
will identify $\Gal(f)$ with its image and consider it as a
permutation group of $\RR_f$. Clearly, $\Gal(f)$ is transitive if
and only if $f$ is irreducible in $K[x]$. Further, we assume that
$p$ does not divide $n$ and $K$ contains a primitive $q$th root of
unity $\zeta$.

We write (as in \cite[\S 3] {ZarhinCrelle}) $V_{f,p}$ for the
$(n-1)$-dimensional $\F_p$-vector space of functions $\{\phi:\RR_f
\to \F_p, \ \sum_{\alpha\in \RR_f}\phi(\alpha)=~0\}$ provided with
a natural {\sl faithful} action of the permutation group
$\Gal(f)\subset \Perm(\RR_f)$. It is  the {\sl heart} over the
field $\F_p$ of the group $\Gal(f)$ acting on the set $\RR_f$
 \cite{Mortimer,ZarhinCrelle}.

\begin{rem}
\label{doubleV}
 It is well-known \cite{Klemm} (see also \cite{ZarhinCrelle,ZarhinMZ})
  that if $\Gal(f)$ acts {\sl doubly
 transitively} on $\RR_f$ then the centralizer
 $$\End_{\Gal(f)}(V_{f,p})=\F_p.$$
 The surjection $\Gal(K)\twoheadrightarrow \Gal(f)$ provides
 $V_{f,p}$ with the natural structure of $\Gal(K)$-module. Clearly,
if $\Gal(f)$ acts {\sl doubly
 transitively} on $\RR_f$ then the centralizer
 $$\End_{\Gal(K)}(V_{f,p})=\F_p.$$
\end{rem}

 Let $C=C_{f,q}$ be the smooth projective model of the smooth
affine $K$-curve
            $y^q=f(x)$.
So $C$ is a smooth projective curve defined over $K$. and the
rational function $x \in K(C)$ defines a finite cover  $\pi:C \to
\P^1$ of degree $q$. Let $\B'\subset C(\bar{K})$ be the set of
ramification points.  Clearly, the restriction of $\pi$ to $\B'$
is an {\sl injective} map $\B' \hookrightarrow \P^1(\bar{K})$,
whose image is the disjoint union of $\infty$ and  $\RR_f$. We
write
$$\B=\pi^{-1}(\RR_f)=\{(\alpha,0)\mid \alpha \in \RR_f\} \subset \B' \subset C(\bar{K}).$$
Clearly, $\pi$ is ramified at each point of $\B$ with ramification
index $q$. The set $\B'$ is the disjoint union of $\B$ and a
single point $\infty':=\pi^{-1}(\infty)$. In addition, the
ramification index of $\pi$ at $\pi^{-1}(\infty)$ is also $q$. By
Hurwitz's formula,  the genus $g=g(C)=g(C_{f,q})$ of $C$ is
$(q-1)(n-1)/2$.

\begin{lem}
\label{genre}
 Let us  consider  the plane triangle
(Newton polygon)
$$\Delta_{n,q}:=\{(j,i)\mid 0\le j,\quad 0\le i, \quad qj+ni\le nq\}$$
with the vertices $(0,0)$, $(0,q)$ and $(n,0)$. Let $L_{n,q}$ be
the set of {\sl integer points} in the interior of $\Delta_{n,q}$.
Then:
\begin{itemize}
\item[(i)]
 If $(j,i)\in L_{n,q}$ then
$$1\le j \le n-1; \quad 1 \le i \le q-1;\quad q(j-1)+(q+1)\le n(q-i).$$
\item[(ii)]
 The genus $g=(q-1)(n-1)/2$ coincides with the number of
elements of $L_{n,q}$. \item[(iii)]
$$\omega_{j,i}:=x^{j-1}dx/y^{q-i}=x^{j-1}y^idx/y^q=x^{j-1}y^{i-1} dx/y^{q-1}$$
is a differential of the first kind on $C$ for each $(j,i) \in
L_{n,q}$.
\end{itemize}
\end{lem}
\begin{proof} (i) is obvious.
In order to prove (ii), notice that $(q-1)(n-1)$ is the number of
interior integer points in the rectangular
$$\Pi_{n,q}:=\{(j,i)\mid 0\le j\le n, 0\le i\le q\},$$
whose diagonal that connects $(0,0)$ and $(n,q)$ does not contain
interior integer points, because $n$ and $q$ are relatively prime.
In order to get the {\sl half} of $q-1)(n-1)$, one should notice
that the map $(j,i) \mapsto (n-j,q-i)$  establishes a bijection
between $L_{n,q}$ and the set of interior integer points of
$\Pi_{n,q}$ outside $L_{n,q}$. Let us prove (iii). Clearly,
$\omega_{j,i}$ has no poles outside $\B'$. For each
$\alpha\in\RR_f$ the orders of zero of $dx =d(x-\alpha)$ and $y$
at $\pi^{-1}(\alpha,0)$ are $q-1$ and $1$ respectively. Since
$(q-1)-i\ge 0$, the differential $\omega_{j,i}$ has no pole at
$\pi^{-1}(\alpha,0)$. The orders of pole of $x,dx,y$ at $\infty$
are $q,q+1,n$ respectively. Since $n(q-i)\ge q(j-1)+(q+1)$, the
differential $\omega_{j,i}$ has no pole at $\infty$. This implies
easily that the collection $\{\omega_{j,i}\}_{(j,i)\in L_{n,q}}$
is a basis in the space of differentials of the first kind on $C$.
\end{proof}

 There is a non-trivial birational
 $\bar{K}$-automorphism of $C$
 $$\delta_q:(x,y) \mapsto (x, \zeta y).$$
Clearly, $\delta_q^q$ is the identity map and
 the set of fixed points of $\delta_q$ coincides with $\B'$.
 Let us consider $C_(f,q)$ as curve over $\bar{K}$ and let
$\Omega^1(C_{(f,q)})$ be the ($g(C_{f,q})$-dimensional) space of
differentials of the first kind on $C$.  By functoriality,
$\delta_q$ induces on $\Omega^1(C_{(f,q)})$  a certain
$\bar{K}$-linear automorphism $\delta_q^*:\Omega^1(C_{(f,q)})\to
\Omega^1(C_{(f,q)})$.

\begin{lem}
 \label{dfk}
 \begin{itemize}
 \item[(i)]
The collection $\{\omega_{j,i}=x^{j-1} dx/y^{q-i}\mid (i,j) \in
L_{n,q}\}$ is an eigenbasis of $\Omega^1(C_{(f,q)})$ with respect
to  $\delta_q^*$; an eigenvector $\omega_{j,i}=x^{j-1} dx/y^{q-i}$
corresponds to eigenvalue $\zeta^{i}$.
\item[(ii)] If a positive
integer $i$ satisfies $i<q$ then $\zeta^{-i}$ is an eigenvalue of
of $\delta_q^*$ if and only if  $\left[\frac{ni}{q}\right]>0$. If
this is the case then the multiplicity of eigenvalue $\zeta^{-i}$
is $\left[\frac{ni}{q}\right]$.
 \item[(iii)] $1$ is not an eigenvalue of $\delta_q^*$.
 \item[(iv)]
${\mathcal P}_q(\delta_q^*)={\delta_q^*}^{q-1}+\cdots
+\delta_q^*+1=0$ in $\End_{\bar{K}}(\Omega^1(C_{f,q)}))$.
 \item[(v)]
 If $q=p^r$ then $\zeta^{-i}$ is an eigenvalue of $\delta_q^*$ for
each integer $i$ with $p^r-p^{r-1}\le i \le p^r-1=q-1$.
\item[(vi)] If ${\mathcal H}(t)$ is a polynomial with rational
coefficients such that ${\mathcal H}(\delta_q^*)=0$ in
$\End_{\bar{K}}(\Omega^1(C_{f,q}))$ then ${\mathcal H}(t)$ is {\sl
divisible} by ${\mathcal P}_q(t)$ in $\Q[t]$.
\end{itemize}
\end{lem}

\begin{proof}
The assertion (i) follows easily from Lemma \ref{genre}. This
implies (iii), which, in turn, implies (iv). The assertion (i)
implies that if $1\le i<q$ then $\zeta^{-i}$ is an eigenvalue of
$\delta_q^*$ if and only if the corresponding (to $q-i$)
horizontal line contains an interior integer point; if this is the
case then the multiplicity coincides with the number of interior
points on this line. Clearly, this number of points is
$\left[\frac{ni}{q}\right]$. In order to prove (v), recall that
$n\ge 3$. This implies that
$$ni\ge 3(p^r-p^{r-1})=3\frac{p-1}{p}q\ge \frac{3}{2}q>q$$ and
therefore  $\left[\frac{ni}{q}\right]\ge 1>0$, which proves (v),
in light of (ii). In order to prove (vi), it suffices to notice
that, thanks to (v), we know that among the eigenvalues of
$\delta_q^*$ there is a primitive $q$th root of unity.
\end{proof}

\begin{rem}
\label{model} In order to describe $C$ explicitly, we use a
construction from \cite[pp. 3358-3359]{Towse} (see also \cite[Ch.
3a, Sect.1]{Mumford2}). Let us put $\tilde{f}(x)=x^n f(1/x)\in
K[x]$ and pick positive integers $a,b$ with $bn-aq=1$. (Such $a,b$
do exist, since $n$ and $q$ are relatively prime.) Let us consider
the birational transformation
$$x=\frac{1}{s^bt^q}, y=\frac{1}{s^a}{t^n}; s=x^{-n}y^{q}=x^{-n}f(x), t=x^{a}y^{-b}. $$
If $f_0\ne 0$ is the leading coefficient of $f(x)$ then
$$W_1:y^q=f(x)=f_0\prod_{\alpha\in\RR}(x-\alpha)$$ becomes
$$W_2:s=\tilde{f}(s^bt^q)=f_0\prod_{\alpha\in\RR}(1-\alpha s^bt^q)).$$
There are no points  with $s=0$ (and therefore for every
$\alpha\in \RR_f$ there are no points  with $\alpha s^bt^q=1$). When
$s\ne 0,t\ne 0$, we are on the original affine piece $W_1$ of the
curve given in terms of $x$ and $y$. (When $t=0$ we get the point
$\infty'=(f_0,0)$.) The automorphism $\delta_q$ sends $(s,t)$ to
$(s,\zeta^{-b}t)$. If we glue  together the affine curves $W_1$
and $W_2$, identifying $\{x\ne 0,y\ne 0\}$ and $\{s\ne 0,t\ne 0\}$
via the birational transformation above, then we get a smooth
$K$-curve $C'$. The finite maps of affine $K$-curves
$$x: W_1\to \A^1, \ s^b t^q: W_2 \to  \A^1 \setminus
\{\alpha^{-1}\mid\alpha\in\RR_f,\alpha\ne 0\}$$ are gluing together
to a finite map $C'\to \P^1$, since $\P^1$ coincides with the
union of the images of the open embeddings
$$\A^1\hookrightarrow\P^1, t\mapsto (t:1),\ \A^1 \setminus
\{\alpha^{-1}\mid\alpha\in\RR_f,\ \alpha\ne 0\}\hookrightarrow\P^1,
t\mapsto (1:t).$$ It follows that $C'$ is proper over $K$ and
therefore is projective (see also \cite[Ch. III, Ex.
5.7]{Hartshorne}). This implies that $C'$ is biregularly
isomorphic to $C$ over $K$.
\end{rem}

Let $J(C_{f,q})=J(C)$ be the jacobian of $C$. It is a
$g$-dimensional abelian variety defined
 over $K$ and one may view (via Albanese functoriality) $\delta_q$ as an element of
 $\Aut(C) \subset\Aut(J(C)) \subset \End(J(C))$
such that
  $\delta_q \ne \I$ but $\delta_q^q=\I$
where $\I$ is the identity endomorphism of $J(C)$. We write
$\Z[\delta_q]$ for the subring of $\End(J(C))$ generated by
$\delta_q$.

\begin{rems}
\label{alb}
\begin{itemize}
 \item[(i)]
 The point $\infty'\in C_{f,p}(K)$ is one of  $\delta_q$-invariant
 points. The map
$$\tau: C_{f,q} \to J(C_{f,q}), \quad P\mapsto \cl((P)-(\infty'))$$ is
an embedding of smooth projective algebraic varieties and it is
well-known \cite[Sect. 2.9]{Shimura} that the induced map $\tau^*:
\Omega^1(J(C_{f,q})) \to \Omega^1(C_{f,q})$ is an isomorphism
obviously commuting with the actions of $\delta_q$. (Here $\cl$
stands for the linear equivalence class.) This implies that
$n_{\sigma_i}$ coincides with the dimension of the eigenspace of
$\Omega^1(C_{(f,q)})$ attached to the eigenvalue $\zeta^{-i}$ of
$\delta_q^*$. Applying Lemma \ref{dfk}, we conclude that if
${\mathcal H}(t)$ is a monic polynomial with integer coefficients
such that ${\mathcal H}(\delta_q)=0$ in $\End(J(C_{f,q}))$ then
${\mathcal H}(t)$ is divisible by ${\mathcal P}_q(t)$ in $\Q[t]$
and therefore in $\Z[t]$.
 \item[(ii)]
The automorphism $\delta_q:C_{f,p}\to C_{f,p}$ induces, by Picard
functoriality, the automorphism $\delta_q{^{\prime}}:J(C_{f,q})\to
J(C_{f,q})$ that sends $\cl((P)-(\infty'))$ to
$\cl(\delta_q^*((P)-(\infty')))=\cl((\delta_q{^{-1}}(P)-(\infty'))$.
This implies that
$$\delta_q{^{\prime}}=\delta_q{^{-1}} \in \Aut(J(C_{f,q}))\subset
\End(J(C_{f,q})).$$
 \end{itemize}
\end{rems}

\begin{rem}
\label{multprim}  Clearly, the set $S$ of eigenvalues $\lambda$ of
$\delta_q^*:\Omega^1(J(C_{f,q})) \to \Omega^1(J(C_{f,q}))$ with
$\PP_{q/p}(\lambda)\ne 0$  consists of {\sl primitive} $q$th roots
of unity $\zeta^{-i}$ ($1 \le i<q, (i,p)=1$) with
$\left[\frac{ni}{q}\right]>0$ and the multiplicity of $\zeta^{-i}$
equals  $\left[\frac{ni}{q}\right]$, thanks to Remarks \ref{alb}
and Lemma \ref{dfk}. Let us consider the sum
$$M=\sum_{1 \le i<q,(i,p)=1}\left[\frac{ni}{q}\right]$$
 of multiplicities of eigenvalues from $S$. Then
$$M=(n-1)\frac{\varphi(q)}{2}=\frac{(n-1)(p-1)p^{r-1}}{2}.$$

See \cite[Remark 4.6 on p. 353]{ZarhinM} for the proof in which
one has only to replace references to \cite[Remarks 4.5 and
4.4]{ZarhinM} by references to Remark \ref{alb}(i) and Lemma
\ref{dfk} respectively.

Clearly, if the abelian (sub)variety
$Z:=\PP_{q/p}(\delta_q)(J(C_{f,q}))$ has dimension $M$ then the
data $Y=J(C_{f,q}), \delta=\delta_q, P=\PP_{q/p}(t)$ satisfy the
conditions of Theorem 3.10 of \cite{ZarhinM}.
\end{rem}

\begin{lem}
\label{order}  Let $D=\sum_{P\in \B}a_P (P)$ be a divisor on
$C=C_{f,p}$ with degree $0$ and support in $\B$. Then $D$ is
principal if and only if all the coefficients $a_P$ are divisible
by $q$.
\end{lem}

\begin{proof}
See \cite[Lemma 4.7 on p. 354]{ZarhinM}.
\end{proof}

\begin{lem}
\label{cycl}
 $1+\delta_q+ \cdots + \delta_q^{q-1}=0$ in $\End(J(C_{f,q}))$.
The subring $\Z[\delta_q] \subset \End(J(C_{f,q}))$ is isomorphic
to the ring $\Z[t]/{\mathcal P}_q(t)\Z[t]$.  The $\Q$-subalgebra
$\Q[\delta_q]\subset\End^0(J(C_{f,q}))=\End^0(J(C_{f,q}))$ is
isomorphic to $\Q[t]/{\mathcal
P}_q(t)\Q[t]=\prod_{i=1}^r\Q(\zeta_{p^i})$.
\end{lem}

\begin{proof}
If $q=p$ is a prime this assertion is proven in
 \cite[p.~149]{Poonen},  \cite[p.~458]{SPoonen}.
 (If $n>3$ then the assertion is proven in \cite{ZarhinM}.) So, further we
 may assume that $q>p$.
 The group $J(C_{f,q}))(\bar{K})$ is generated by
 divisor classes of the form $(P)-(\infty)$ where $P$ is a finite
 point on $C_{f,p}$. The divisor of the rational function $x-x(P)$ is
 $(\delta_q^{q-1}P)+\cdots + (\delta_q P) + (P)- q(\infty)$. This
 implies that
 ${\mathcal P}_q(\delta_q)=0 \in \End(J(C_{f,q}))$.
 Applying Remark \ref{alb}(ii), we conclude that ${\mathcal P}_q(t)$
 is the minimal polynomial of $\delta_q$ in $\End(J(C_{f,q}))$.
\end{proof}

Let us define the abelian (sub)variety
$$J^{(f,q)}:= \PP_{q/p}(\delta_q)(J(C_{f,q}))\subset J(C_{f,q}).$$
Clearly, $J^{(f,q)}$ is a $\delta_q$-invariant abelian subvariety
defined over $K$. In addition, $\Phi_q(\delta_q)(J^{(f,q)})=0$.

\begin{rem}
\label{qp} If $q=p$ then $\PP_{q/p}(t)=\PP_{1}(t)=1$ and therefore
$J^{(f,p)}=J(C_{f,p})$.
\end{rem}

\begin{rem}
\label{nonzero} Since the polynomials $\Phi_q$ and $\PP_{q/p}$ are
relatively prime, the homomorphism $\PP_{q/p}(\delta_q):J^{(f,q)}
\to J^{(f,q)}$  has finite kernel and therefore is an isogeny. In
particular, it is surjective.
\end{rem}

\begin{lem}
\label{fixP}
 $\dim (J^{(f,q)})=\frac{(p^r-p^{r-1})(n-1)}{2}$
and there is an $K$-isogeny
  $J(C_{f,q})\to
J(C_{f,q/p})\times J^{(f,q)}$. In addition, the Galois modules
$V_{f,p}$ and $$(J^{(f,q)})^{\delta_q}:=\{z \in
J^{(f,q)}(\bar{K})\mid \delta_q(z)=z\}$$ are isomorphic.
\end{lem}

\begin{proof}
Let us consider the curve $C_{f,q/p}:y_1^{q/p}=f(x_1)$ and a
regular surjective map $\pi_1:C_{f,q} \to C_{f,q/p}, \quad x_1=x,
y_1=y^p$. Clearly, $\pi_1\delta_q=\delta_{q/p}\pi_1$. By Albanese
functoriality, $\pi_1$ induces a certain surjective homomorphism
of jacobians $J(C_{f,q}) \twoheadrightarrow J(C_{f,q/p})$ which we
continue to denote by $\pi_1$. Clearly, the equality
$\pi_1\delta_q=\delta_{q/p}\pi_1$ remains true in
$\Hom(J(C_{f,q}),J(C_{f,q/p}))$. By Lemma \ref{cycl}, $${\mathcal
P}_{q/p}(\delta_{q/p})=0 \in \End(J(C_{f,{q/p}})).$$ It follows
from Lemma \ref{nonzero} that $\pi_1(J^{(f,q)})=0$ and therefore
$\dim(J^{(f,q)})$ does not exceed
$$\dim(J(C_{f,q}))-\dim(J(C_{f,q/p}))=$$
$$\frac{(p^r-1)(n-1)}{2}-\frac{(p^{r-1}-1)(n-1)}{2}=\frac{(p^r-p^{r-1})(n-1)}{2}.$$
By definition of $J^{(f,q)}$, for each divisor $D=\sum_{P\in
\B}a_P (P)$ the linear equivalence class of $\$p^{r-1}D=\sum_{P\in
\B}p^{r-1}a_P(P)$ lies in  $(J^{(f,q)})^{\delta_q}\subset
J^{(f,q)}(\bar{K})\subset J(C_{f,q})(\bar{K})$. It follows from
Lemma \ref{order} that the class of $p^{r-1}D$ is zero if and only
if all $p^{r-1}a_P$ are divisible by $q=p^r$, i.e. all $a_P$ are
divisible by $p$. This implies that the set of linear equivalence
classes of $p^{r-1}D$ is a Galois submodule isomorphic to
$V_{f,p}$. We want to prove that $(J^{(f,q)})^{\delta_q}=V_{f,p}$.

Recall that $J^{(f,q)}$ is $\delta_q$-invariant and the
restriction of $\delta_q$ to $J^{(f,q)}$ satisfies the $q$th
cyclotomic polynomial. This allows us to define the homomorphism
$\Z[\zeta_q] \to \End(J^{(f,q)})$ that sends $1$ to the identity
map and $\zeta_q$ to $\delta_q$. Let us put $$E=\Q(\zeta_q),
\O=\Z[\zeta_q]\subset \Q(\zeta_q)=E.$$ It is well-known that $\O$
is the ring of integers in $E$, the ideal
$\lambda=(1-\zeta_q)\Z[\zeta_q]=(1-\zeta_q)\O$ is  maximal in $\O$
with $\O/\lambda=\F_p$ and $\O\otimes\Z_p=\Z_p[\zeta_q]$ is the
ring of integers in the field $\Q_p(\zeta_q)$. Notice also that
$\O\otimes\Z_p$ coincides with the completion $\O_{\lambda}$ of
$\O$ with respect to the $\lambda$-adic topology and
$\O_{\lambda}/\lambda \O_{\lambda}=\O/\lambda=\F_p$.

 It follows  from Lemma 3.3  of \cite{ZarhinM} that
$$d=\frac{2\dim(J^{(f,q)})}{[E:\Q]}=\frac{2\dim(J^{(f,q)})}{p^r-p^{r-1}}$$
is a positive integer, the $\Z_p$-Tate module $T_p(J^{(f,q)})$ is
a free $\O_{\lambda}$-module of rank $d$. (See also \cite[Prop.
2.2.1 on p, 769]{Ribet}.) Using the displayed formula (5) on p.
347 in  \cite[\S 3]{ZarhinM},
 we conclude that
$$(J^{(f,q)})^{\delta_q}=\{u \in J^{(f,q)}(\bar{K})\mid
(1-\delta_q)(u)=0\}=J^{f,q}_{\lambda}=T_p(J^{f,q})\otimes_{\O_{\lambda}}\F_p$$
is a $d$-dimensional $\F_p$-vector space.
 Since $(J^{(f,q)})^{\delta_q}$ contains $(n-1)$-dimensional
$\F_p$-vector space $V_{f,p}$, we have $d \ge n-1$. This implies
that $$2\dim(J^{(f,q)})=d (p^r-p^{r-1}) \ge (n-1)(p^r-p^{r-1})$$
and therefore $$\dim(J^{(f,q)}) \ge
\frac{(n-1)(p^r-p^{r-1})}{2}.$$ But we have already seen that
$$\dim(J^{(f,q)}) \le \frac{(n-1)(p^r-p^{r-1})}{2}.$$ This implies
that $$\dim(J^{(f,q)}) = \frac{(n-1)(p^r-p^{r-1})}{2}.$$ It
follows that $d=n-1$ and therefore
$(J^{(f,q)})^{\delta_q}=V_{f,p}$. Dimension arguments  imply that
$J^{(f,q)}$ coincides with the identity component of $\ker(\pi_1)$
and  therefore there is a $K$-isogeny between
  $J(C_{f,q})$ and $J(C_{f,q/p})\times J^{(f,q)}$.
\end{proof}

\begin{cor}
\label{split} There is a $K$-isogeny $$J(C_{f,q})\to
J(C_{f,p})\times \prod_{i=2}^r J^{(f,p^i)}=\prod_{i=1}^r
J^{(f,p^i)}.$$
\end{cor}
\begin{proof}
Combine Corollary \ref{fixP}(ii) and Remark \ref{qp} with easy
induction on $r$.
\end{proof}

\begin{thm}
\label{maximalV}
 Suppose that $n \ge 3$ is an integer. Let $p$ be a
prime, $r \ge 1$ an integer and $q=p^r$.  Let $f(x)\in K[x]$ be a
polynomial of degree $n$ without multiple roots.  Then the image
$\O$ of $\Z[\delta_q] \to \End(J^{(f,q)})$ is isomorphic to
$\Z[\zeta_q]$ and the Galois module
$(J^{(f,q)})^{\delta_q}=J^{(f,q)}_{\lambda}$ is isomorphic to
$V_{f,p}$. If $\Gal(f)$ acts {\sl doubly
 transitively} on $\RR_f$ then the centralizer
 $$\End_{\Gal(K)}(J^{(f,q)}_{\lambda})=\F_p.$$
\end{thm}

\begin{proof}
 Clearly, $\O$ is isomorphic to $\Z[\zeta_q]$. Let us
put $\lambda=(1-\zeta_q)\Z[\zeta_q]$. By Lemma \ref{fixP}(iii),
the Galois module $(J^{(f,q)})^{\delta_q}=J^{(f,q)}_{\lambda}$ is
isomorphic to $V_{f,p}$. The last assertion follows from Remark
\ref{doubleV}.
\end{proof}

\begin{rem}
\label{multJ} It follows from Theorem \ref{maximalV} that there is
an embedding
$$\Z[\zeta_q] \hookrightarrow \End_K(J^{(f,q)})\subset \End(J^{(f,q)}), \ \zeta_q \mapsto \delta_q:
J^{(f,q)}\to J^{(f,q)},$$ which sends $1$ to the identity
automorphism of $J^{(f,q)}$. This embedding extends by
$\Q$-linearity to the embedding
$$\Q(\zeta_q) \hookrightarrow  \End^0(J^{(f,q)}).$$
 On the other hand, let us consider the
induced linear operator $$\delta_q^*:\Omega^1(J^{(f,q)})\to
\Omega^1(J^{(f,q)}).$$ It follows from Theorem 3.10 of
\cite{ZarhinM} combined with Remark \ref{multprim} that its
spectrum consists of primitive $q$th roots of unity $\zeta^{-i}$
($1\le i<q$, the prime $p$ does not divide $i$) with $[ni/q]>0$
and the multiplicity of $\zeta^{-i}$ equals $[ni/q]$. It follows
that if $\sigma_i:E=\Q(\zeta_q)\hookrightarrow \C$ is the field
embedding that sends $\zeta_q$ to $\zeta^{-i}$  and
$$\Omega^1(J^{(f,q)})_{\sigma_i}=\{x \in
\Omega^1(J^{(f,q)})\mid ex=\sigma_i(e)x \quad \forall e\in E\}$$
then
$$\dim_{\bar{K}}(\Omega^1(J^{(f,q)})_{\sigma_i})
=\left[\frac{ni}{q}\right].$$
\end{rem}

\begin{lem}
\label{idcomp} The abelian subvariety $J^((f,q))\subset
J(C_{f,q})$ coincides with the identity component of the kernel of
$$\Phi_q(\delta_q): J(C_{f,q})\to J(C_{f,q}).$$
\end{lem}

\begin{proof}
We know that $J^{(f,q)}\subset\ker(\Phi_q(\delta_q)):=W$. Let
$\pi_1:J(C_{f,q})\to J(C_{f,q/p})$ be as in the proof of Lemma
\ref{fixP}. Clearly, $J^{(f,q)}$  is the identity component of
$\ker(\pi_1)$ and the image  $\pi_1(W)\subset J(C_{f,q/p})$ is
killed by both $\Phi_q(\delta_q)$ and $\PP_{q/p}(\delta_q)$. Since
the polynomials $\Phi_q$ and $\PP_{q/p}$ are relatively prime,
$\pi_1(W)$ is killed by a non-zero integer (that is their
resultant) and therefore is finite. It follows that $W/J^{(f,q)}$
is also finite and therefore $J^{(f,q)}$ is the identity component
of $W=\ker(\Phi_q(\delta_q))$.
\end{proof}

\begin{rem}
\label{idcompalb} Recall (Remark \ref{alb}) that
$\delta_q^{\prime}=\delta_q^{-1}$ and $\delta_q^q$ is the identity
automorphism of $J(C_{f,q})$. It follows from the displayed
formula (4) at the beginning of this Section that
$\Phi_q(\delta_q)=\Phi_q(\delta_q^{\prime})$ in
$\End(J(C_{f,q}))$. Applying Lemma \ref{idcomp}, we conclude that
$J^{(f,q)}$ is the identity component of
$\ker(\Phi_q(\delta_q^{\prime}))$.
\end{rem}

\begin{thm}
\label{tq4n3}

Suppose that $n=3, q=4$. Then:
\begin{itemize}
\item[(i)]
 $J^{(f,4)}$ is isomorphic to a product of two mutually isogenous
elliptic curves with complex multiplication by $\Q(\sqrt{-1})$.
\item[(ii)] $J^{(f,4)}$ is isogenous to the square of the elliptic
curve $y^2=x^3-x$. \item[(iii)]
 $\End^0(J^{(f,4)})\cong\Mat_2(\Q(\sqrt{-1}))$.
 \end{itemize}
\end{thm}

\begin{proof}
The assertion follows readily from Lemma \ref{q4n3} combined with
Remark \ref{multJ}.

\end{proof}
\begin{thm}
\label{Jfqdouble} Suppose that $k$ is an algebraically closed
field of characteristic zero, $K=k(t)$ and $(n,q) \ne (3,4)$.
Suppose that $\Gal(f)$ acts doubly transitively on $\RR_f$.

Then:

\begin{itemize}
\item[(i)]
$J^{(f,q)}$ is completely non-isotrivial.

\item[(ii)] If $n=3$ or $4$ then $\End(J^{(f,q)})=\Z[\zeta_q]$. In
particular, $J^{(f,q)}$ is absolutely simple and
$\End^0(J^{(f,q)})=\Q(\zeta_q)$.
\end{itemize}
\end{thm}

\begin{proof}
Let us fix a primitive $q$th root of unity $\zeta\in k$ and
consider the subfield $k_f\subset k$ that is obtained by adjoining
to $\Q$ all coefficients of $f(x)$ and $\zeta$. We have
$$f(x)\subset k_f(t)[x]\subset k(t)[x]=K[x], \ k_f(t)(\RR_f)\subset K(\RR_f)\subset \bar{K}.$$
Let $k_0$ be the algebraic closure of $k_f$ in $k_f(t)(\RR_f)$.
Since a subfield of a finitely generated overfield is also
finitely generated \cite[Th. 4.3.6]{Rotman}, $k_0$ is finitely
generated over $k_f$ and therefore is a finite algebraic extension
of $k_f$ (see \cite[Chap. V, Sect. 14, no. 7, Cor. 1]{B}).
Clearly, $k_0(t)(\RR_f)=\ k_f(t)(\RR_f)$; in particular, $k_0$ is
algebraically closed in $k_0(t)(\RR_f)$. Since $k_0$ lies in
$\bar{K}=\overline{k(t)}$ and is algebraic over $k_f\subset k$, we
have $k_0\subset k$. We have
$$\Q \subset k_f \subset k_0\subset k, \ f(x)\in k_f(t)[x]\subset k_0(t)[x]
\subset k(t)[x];$$ in addition, $k_0(t)(\RR_f)$ is the {\sl
splitting field} of $f(x)$ over $k_0(t)$. Clearly, $k_0$ is
finitely generated over $\Q$ and contains $\zeta$.

Let $\bar{k}_0$ be the algebraic closure of $k_0$ in $k$. Since
$k$ is algebraically closed, $\bar{k}_0$ is also algebraically
closed and therefore $\bar{k}_0(t)$ is algebraically closed in
$k(t)$. Since the field extension
$\bar{k}_0(t)(\RR_f)/\bar{k}_0(t)$ is algebraic,
$\bar{k}_0(t)(\RR_f)$ and $k(t)$ are linearly disjoint over
$\bar{k}_0(t)$. It follows that $f(x)$ has the same Galois group
$\Gal(f)$ over $k(t)$ and $\bar{k}_0(t)$.

Since $k_0$ is algebraically closed in $k_0(t)(\RR_f)$, the fields
$k_0(t)(\RR_f)$ and $\bar{k}_0$ are linearly disjoint over $k_0$.
This implies that $f(x)$ has the same Galois group over $k_0(t)$
and $\bar{k}_0(t)$. It follows that the Galois group of $f(x)$
over  $k_0(t)$ coincides with $\Gal(f)$; in particular, it acts
doubly transitively on $\RR_f$.

 Since $k_0$ is finitely generated over
$\Q$, it is a {\sl hilbertian} field \cite[pp. 129--130]{SerreM}.
Let us put $K_0=k_0(t)$. An elementary substitution allows us to
assume that $f(x)=:F(t,x)\in k_0[t][x]=k_0[t,x]$. (Indeed, if
$f(x)=G(x,t)/h(t)$ with $G(x,t)\in k_0[t,x], \ h(t)\in k_0[t]$
then $C_{f,q}$ is $K_0$-birationally isomorphic to the
superelliptic curve $w^q=h(t)^{q-1}G(x,t)$ where $w=h(t)\cdot y$.)

Using Remark \ref{model}, one may easily check that there exists a
finite set $B$ of closed points of the affine $k_0$-line
$\A^1_{k_0}$ and a smooth proper $\A^1_{k_0}\setminus B$-curve
$c:\CV\to \A^1_{k_0}\setminus B$ that enjoy the following
properties:

\begin{itemize}
\item For all $a\in k\setminus B$ the polynomial $f_a(x):=F(a,x)$
has no multiple roots; \item The ``affine part" of  $c:\CV\to
\A^1_{k_0}\setminus B$ is defined by equation $y^q=F(x,t)$;
 \item
The generic fiber of  $c$ coincides with $C_{f,q}$ and for all
$a\in k\setminus B$ the corresponding closed fiber is $C_{f_a,q}$.
In addition, the infinite points $\infty'$ give rise to a section
$$\boldsymbol{\infty}':\A^1_{k_0}\setminus B \to \CV.$$
\item
 The map $(x,y)\mapsto (x,\zeta y)$ gives rise to a periodic
automorphism $\boldsymbol{\delta}_q:\CV\to \CV$ of the
$\A^1_{k_0}\setminus B$-scheme that coincides with $\delta_q$ on
the generic fiber and all closed fibers of $c$; \item
$\boldsymbol{\delta}_q$ acts as the identity map on the image of
$\boldsymbol{\infty}'$.
\end{itemize}

 Taking the (relative) jacobian of $c:\CV\to \A^1_{k_0}\setminus B$ \cite[Ch. 9]{Neron}, we get an
abelian scheme $\J \to \A^1_{k_0}\setminus B$, whose generic fiber
coincides with $J(C_{f,q})$ and for each $a\in k_0\setminus B$ the
corresponding closed fiber is $J(C_{f_{a},q})$.

Since $k_0$ is hilbertian and algebraically closed in
$k_0(t)(\RR_f)$, there exist two distinct $a,d\in k_0\setminus B$
such that  both $\Gal(f_{a})$ and $\Gal(f_{d})$ coincide with
$\Gal(f)$ and the splitting fields of $f_{a}$ and $f_{d}$ are
linearly disjoint over $k_0$ \cite[Sect. 10.1]{SerreM}. In
particular, both $\Gal(f_{a})$ and $\Gal(f_{d})$ are doubly
transitive. It follows from Remark \ref{doubleV} that the
$\Gal(k_0)$-modules $V_{f_{a}}$ and $V_{f_{d}}$ satisfy
$$\End_{\Gal(k_0)}(V_{f_{a}})=\F_p, \ \End_{\Gal(k_0)}(V_{f_{d}})=\F_p \eqno (5).$$

First, assume that $q=p$. It follows from Corollary 3.6 in
\cite{ZarhinMZ} applied to $X=J(C_{f_{a},p}), Y=J(C_{f_{d},p}),
E=\Q(\zeta_p), \lambda=(1-\zeta_p)$ that either there are no
non-zero $\bar{k}_0$-homomorphisms between $J(C_{f_{a},p})$ and
$J(C_{f_{d},p})$ or the centralizer of $\Q(\zeta_p)$ in
$\End^0(J(C_{f_{a},p}))$ has $\Q(\zeta_p)$-dimension
$$\left(\frac{2\dim(J(C_{f_{a},p}))}{[\Q(\zeta_p):\Q]}\right)^2=(n-1)^2.$$
Combining  Theorem \ref{cyclmult}(b) with Lemma \ref{dfk}, we
conclude that all $\bar{k}_0$-homomorphisms between
$J(C_{f_{a},p})$ and $J(C_{f_{d},p})$ are zero. Since $\bar{k}_0$
is algebraically closed, all $k$-homomorphisms between
$J(C_{f_{a},p})$ and $J(C_{f_{d},p})$ are also zero.

It follows from Lemma \ref{nonconst} applied to
$X=J(C_{f,p})=J^{(f,p)}$ and $\X=\J \times_{{k}_0} k$ that
$J(C_{f,p})$ is completely non-isotrivial. This proves (i) when
$q=p$.

 Suppose that $n=3$ or $4$. Since every
abelian variety of CM-type is defined over $\bar{\Q}$ \cite{Oort},
Theorem \ref{bigend} combined with the complete non-isotriviality
of $J^{(f,p)}$ imply that
 $\Q(\zeta_p)$ coincides with its own centralizer in
$\End^0(J(C_{f,p}))$. Then $\Q(\zeta_p)$ contains the center of
$\End^0(J(C_{f,p}))$. It follows from Theorem \ref{cyclmult}(a)
that $\Q(\zeta_p)$ is the center of $\End^0(J(C_{f,p}))$. Then the
centralizer of $\Q(\zeta_p)$ coincides with the whole
$\End^0(J(C_{f,p}))$ and therefore
$\Q(\zeta_p)=\End^0(J(C_{f,p}))$. Since $\Z(\zeta_p]\subset
\End(J(C_{f,p}))$ is the maximal order in $\Q(\zeta_p)$, we
conclude that $\Z[\zeta_p]= \End(J(C_{f,p}))$. This proves (ii)
when $q=p$.

Now, in order to do the case of arbitrary $q$, we need to
construct an abelian scheme over $\A^1_{k_0}\setminus B$,  whose
generic fiber coincides with $J^{(f,q)}$ and for each $a\in
k_0\setminus B$ the corresponding closed fiber is $J^{(f_a,q)}$.
First, since $J^{(f,q)}$ is an abelian subvariety of $J(C_{f,q})$,
it also has good reduction everywhere at $\A^1\setminus B$
\cite[Lemma 2 on p. 182]{Neron}. Let $\J^{(f,q)}$ be  the
schematic closure of $J^{(f,q)}$ in $\J$ \cite[p. 55]{Neron}. It
follows from \cite[Sect. 7.1, p. 175, Cor. 6]{Neron} that
$\J^{(f,q)}$ is a N\'eron model of $J^{(f,q)}$ over
$\A^1_{k_0}\setminus B$.

The automorphism $\boldsymbol{\delta}_{q}:\CV\to \CV$ induces, by
Picard functoriality  the automorphism of the abelian scheme $\J$
and we denote it by $\boldsymbol{\delta}^{\prime}_{q}:\J\to\J$.
Let $\H\subset\J$ be the kernel of
$\Phi_q(\boldsymbol{\delta}^{\prime}_{q})$; it is a closed group
subscheme of $\J$ that is not necessarily an abelian subscheme. We
write $\H_K$ for the generic fiber of $\H$; it is a closed (not
necessarily connected) algebraic $K$-subgroup of $J(C_{f,q})$.

By Remark \ref{idcompalb}, $J^{(f,q)}\subset
\ker(\Phi_q(\delta^{\prime}_{q}))$. This implies that
$J^{(f,q)}\subset \H_K$ and therefore $\J^{(f,q)}\subset\H$. For
all $a\in k_0\setminus B$, the closed fiber $\H_a$ of $\H$
coincides with $\ker(\Phi_q(\delta^{\prime}_{q}))\subset
J(C_{f_{a},q})$; this implies that $\H_a\subset J(C_{f_{a},q})$ is
a closed algebraic subgroup, whose identity component coincides
with $J^{(f_a,q)}$. On the other hand, the closed fiber
${\J^{(f,q)}}_a$ of $\J^{(f,q)}$ at $a$ is an abelian subvariety
in $\H_a$ and has the same dimension as $J^{(f,q)}$, i.e., the
same dimension as $J^{(f_a,q)}$. It follows that
${\J^{(f,q)}}_a=J^{(f_a,q)}$, i.e., $J^{(f_a,q)}$ coincides with
the reduction of $\J^{(f,q)}$ at $a$.

Now one may carry out the same arguments as in the case of $q=p$.
Namely, pick two distinct $a,d\in k_0\setminus B$ such that both
$\Gal(f_{a})$  and $\Gal(f_{d})$  coincide with $\Gal(f)$ and the
splitting fields of $f_{a}$ and $f_{d}$ are linearly disjoint over
$k_0$. Using the displayed formula (5) and applying Corollary 3.6
of \cite{ZarhinMZ} to $X=J^{(f_{a},q)}, Y=J^{(f_{d},q)},
E=\Q(\zeta_q), \lambda=(1-\zeta_q)$, we conclude that either there
are no non-zero $\bar{k}_0$-homomorphisms between $J^{(f_{a},q)}$
and $J^{(f_{d},q)}$ or the centralizer of $\Q(\zeta_q)$ in
$\End^0(J^{(f_{a},q)})$ has $\Q(\zeta_q)$-dimension
$$\left(\frac{2\dim(J^{(f_{a},q)})}{[\Q(\zeta_q):\Q]}\right)^2=(n-1)^2.$$

Combining  Theorem \ref{cyclmult}(b) with Remark \ref{multJ}, we
conclude that there are no non-zero $\bar{k}_0$-homomorphisms
between $J^{(f_{a},q)}$ and $J^{(f_d,q)}$. Again this implies that
there are no non-zero $k$-homomorphisms between $J^{(f_{a},q)}$
and $J^{(f_d,q)}$. It follows from Lemma \ref{nonconst} applied to
$X=J^{(f,q)}$ and $\X=\J^{(f,q)} \times_{{k}_0}k$ that $J^{(f,q)}$
is completely non-isotrivial. This proves (i) for arbitrary $q$.

Suppose that $n=3$ or $4$. As above, Theorem \ref{bigend} combined
with the complete non-isotriviality of $J^{(f,q)}$ imply that
 $\Q(\zeta_q)$ coincides with its own centralizer in
$\End^0(J^{(f,q)})$. Then $\Q(\zeta_q)$ contains the center of
$\End^0(J^{(f,q)})$. It follows from Theorem \ref{cyclmult}(a)
that $\Q(\zeta_q)$ is the center of $\End^0(J^{(f,q)})$. Then the
centralizer of $\Q(\zeta_q)$ coincides with the whole
$\End^0(J^{(f,q)})$ and therefore $\Q(\zeta_q)=\End^0(J^{(f,q)})$.
Since $\Z[\zeta_q]\subset \End(J^{(f,q)})$ is the maximal order in
$\Q(\zeta_q)$, we conclude that $\Z(\zeta_q]= \End(J^{(f,q)})$.
This proves (ii) for arbitrary $q$.
\end{proof}

\begin{proof}[Proof of main results]
Theorems \ref{primep} and \ref{double} follow easily from Theorem
\ref{Jfqdouble}  combined with Corollary \ref{split}. Theorem
\ref{oddp} follows easily from Theorem \ref{Jfqdouble}  combined
with Corollary \ref{split} and Remark \ref{hom}.
\end{proof}

{\sl Proof of Theorem \ref{deg3}.} It follows from Theorem
\ref{Jfqdouble} applied to $(n,q)=(3,2)$  that the elliptic curve
$C_{f,2}=J(C_{f,2})$ has transcendental $j$-invariant and
therefore $\End^0(J(C_{f,2}))=\Q$. (Alternatively, one may use
Lemma \ref{j} instead of Theorem \ref{Jfqdouble}.)
 Combining this with
Lemma \ref{fixP} and Theorem \ref{tq4n3}, we obtain the first
assertion of Theorem \ref{deg3}. The second assertion follows from
the first one combined with Theorem \ref{Jfqdouble}, Lemma
\ref{fixP} and Remark \ref{hom}.

\end{document}